\documentclass[graybox]{svmult}

\usepackage{makeidx}         
\usepackage{graphicx}        
\usepackage{multicol}        
\usepackage[bottom]{footmisc}

\usepackage{newtxtext}       %
\usepackage[varvw]{newtxmath}       
\usepackage{bbm}
\usepackage{xcolor}

\makeindex             

\begin{document}

\title*{Preconditioning of GMRES for Helmholtz problems with quasimodes}
\author{Victorita Dolean\orcidID{0000-0002-5885-1903}, 
Pierre Marchand\orcidID{0000-0002-2522-6837}, 
Axel Modave\orcidID{0000-0002-9145-6585} and 
Timothée Raynaud\orcidID{0009-0007-1780-4135}
}
\institute{Victorita Dolean \at Department of Mathematics and Computer Science, Eindhoven University of Technology, P.O. Box 513, 5600 MB Eindhoven, The Netherlands, \email{v.dolean.maini@tue.nl}
\and Pierre Marchand, Axel Modave and Timothée Raynaud \at POEMS, CNRS, INRIA, ENSTA, Institut Polytechnique de Paris, 91120 Palaiseau, France, \email{pierre.marchand@inria.fr, axel.modave@ensta.fr, and timothee.raynaud@ensta.fr}
}
%
%
\maketitle

\section{Introduction}
\label{sec:intro}

Time-harmonic scalar wave propagation can be modeled by the Helmholtz equation:\vspace{-0.3cm}
\begin{eqnarray*}
    -\Delta u - k^2 u = f,
\end{eqnarray*}
where \( k > 0 \) is the wave number, \(u\) the unknown field, and \(f\) a given source.
Boundary conditions (and/or radiation conditions if the problem is defined in an unbounded domain) are imposed to ensure a well-posed problem.
For numerical solutions, the computational domain is truncated, and the finite element method gives a linear system\vspace{-0.1cm}
\begin{eqnarray}\label{eq:linsys}
    \mathbf{A} \mathbf{u} = \mathbf{b},
\end{eqnarray}
where \( \mathbf{A} \in \mathbb{C}^{N \times N} \) is sparse, non-singular, and typically non-Hermitian.
The unknown and the given source are \( \mathbf{u},\mathbf{b} \in \mathbb{C}^{N} \).
To guarantee the accuracy of the solution and counteract the pollution effect, the system can be very large.

Krylov solvers such as GMRES are well-suited for solving these large non-Hermitian linear systems.
However, the convergence can be extremely slow due to the highly indefinite nature of the problem.
Then, well-designed acceleration techniques, such as domain decomposition methods, are often necessary to ensure convergence within a reasonable time frame.

The numerical solution is even more difficult for quasiresonant Helmholtz problems.
Quasimodes occur when a Helmholtz problem with Sommerfeld radiation conditions is solved in an unbounded domain.
They refer to an increasing sequence of wave numbers such that the norm of the inverse operator increases quickly~\cite{MarchandGalkowskiEtAl2022AGH}:
\begin{definition}[{Quasimodes~\cite[Definition 1.1]{MarchandGalkowskiEtAl2022AGH}}]\label{def:quasimodes}
    Let \(H^1_{\mathrm{loc}}(\Omega)\) be the space of functions locally in \(H^1(\Omega)\), where \(\Omega \subset \mathbb{R}^2\) is unbounded.
    A family of Dirichlet (or Neumann) quasimodes with quality $\epsilon(k)$ is a sequence
    \(
        {\left\{\left( u_j, k_j \right)\right\}}_{j=1}^\infty \subset H^1_{\text{loc}}(\Omega) \times \mathbb{R},
    \)
    where $u_j = 0$ (or $\partial_{\boldsymbol{n}} u_j = 0$) on $\Gamma$, such that $k_j \to \infty$ as $j \to \infty$, and there exists a compact set $K \subset \Omega$ such that, for all $j$,\vspace{-0.1cm}
    \[
        \textnormal{supp}(u_j) \subset K, \quad \|{-\Delta u_j - k_j^2 u_j}\|_{L^2(\Omega)}  \leq \epsilon(k_j), \quad \text{and} \quad \|{u_j}\|_{L^2(\Omega)} = 1.
    \]
\end{definition}
The effect of quasimodes on the GMRES convergence has been studied~\cite{DoleanMarchandEtAl2025CAG,MarchandGalkowskiEtAl2022AGH}, but not in the context of domain decomposition methods.
In this work, we investigate the impact of quasimodes on the GMRES convergence, within the framework of domain decomposition.

\section{GMRES convergence}
\label{sec:gmres}

Here, we briefly present the Generalized Minimal Residual (GMRES) method~\cite{SaadSchultz1986GGM} and its properties, and we introduce a residual estimate based on the harmonic Ritz (HR) values.

GMRES is an iterative algorithm for solving general linear systems such as~\eqref{eq:linsys}, which seeks approximate solutions in the Krylov subspaces\vspace{-0.1cm}
\begin{equation*}
    \mathcal{K}_l(\mathbf{A}, \mathbf{r}_0) := \mathrm{span}\{\mathbf{r}_0, \mathbf{A}\mathbf{r}_0, \mathbf{A}^2\mathbf{r}_0, \ldots, \mathbf{A}^{l-1}\mathbf{r}_0\},
    \quad \text{for $l<N$},
\end{equation*}
where $\mathbf{r}_0 := \mathbf{b} - \mathbf{A}\mathbf{x}_0$ is the initial residual, with $\mathbf{x}_0$ an initial guess.
At each iteration $l < N$, the approximate solution is computed such that the two-norm of the residual \(\mathbf{r}_l := \mathbf{b} - \mathbf{A}\mathbf{x}_l\) is minimized:\vspace{-0.1cm}
\begin{eqnarray}\label{eq:gmres}
    \mathbf{x}_l := \mathrm{argmin}_{\mathbf{x} \in \mathbf{x}_0 + \mathcal{K}_l(\mathbf{A}, \mathbf{r}_0)} \|\mathbf{b} - \mathbf{A}\mathbf{x}\|_2 \in \mathbf{x}_0 + \mathcal{K}_l(\mathbf{A}, \mathbf{r}_0).
\end{eqnarray}

To study the GMRES convergence, it is usual to consider an equivalent problem for~\eqref{eq:gmres} by expressing the residual at iteration $l$ such that:\vspace{-0.1cm}
\begin{equation}
    \label{eq:minResidual}\textstyle
    \|\mathbf{r}_l\|_2 = \min_{q_l \in \mathcal{P}_l^1} \|q_l(\mathbf{A}) \mathbf{r}_0\|_2,
\end{equation}
where $\mathcal{P}_l^1$ is the set of polynomials of degree at most~$l$ with the constraint $q_l(0) = 1$~\cite{Embree2022HDA}. Here, the notation $p_l \in \mathcal{P}_l^1$ refers to the polynomial such that $\mathbf{r}_l = p_l(\mathbf{A}) \mathbf{r}_0$.

Many GMRES convergence bounds are established in the literature, see e.g.~\cite{Embree2022HDA}.
However, most of them fail to capture the non-linear behavior of the convergence rate, which is typically observed for quasiresonant Helmholtz problems~\cite{DoleanMarchandEtAl2025CAG}.
Harmonic Ritz (HR) values, which are eigenvalues estimates of $\mathbf{A}$ over the GMRES iterations~\cite[Chapter 26]{TrefethenEmbree2005SPB}, can be used to interpret the non-linear behavior of the GMRES convergence rate~\cite{Cao1997NCB,DoleanMarchandEtAl2025CAG}.
\begin{definition}[Harmonic Ritz (HR) values]\label{def:harmonicRitzValues}
    At iteration $l > 0$, the \emph{harmonic Ritz values} ${\{\nu_j^{(l)}\}}_{j=1}^l$ are the roots of the minimizing polynomial $p_l \in \mathcal{P}_l^1$ from~\eqref{eq:minResidual}.
\end{definition}

\begin{theorem}[See {\cite[Theorem 2.6]{DoleanMarchandEtAl2025CAG}}]\label{thm:harmonicRitzValues} Let \(\mathbf{A}\) be a non-singular and diagonalizable matrix.
     Let $J \in [1, l]$. Consider:\vspace{-0.15cm}
     \begin{itemize}
        \item a set $N_J^{(l)}$ of $J$ HR values at iteration $l$;
        \item a set $\Lambda_J$ of $J$ eigenvalues of $\mathbf{A}$;
        \item the set $\Lambda_J^c := \sigma(\mathbf{A}) \setminus \Lambda_J$ of the other eigenvalues of $\mathbf{A}$;
        \item \(s_J^{l}(z) :=  {\prod_{\lambda_j \in \Lambda_J} \left(1 - {z}/{\lambda_j}\right)} \cdot {\prod_{\nu_j^{(l)} \in N_J^{(l)}} {\left(1 - {z}/{\nu_j^{(l)}}\right)}^{-1}}\).
     \end{itemize}
    Then, for any $m > 0$ such that $l + m < N$, the local relative residual satisfies\vspace{-0.1cm}
    \begin{eqnarray*}
        \frac{\|\mathbf{r}_{l+m} \|_2}{\|\mathbf{r}_l \|_2}  \leq \left( \sum_{\lambda_i\in \Lambda_J^c} \kappa(\lambda_i) \right) \ \max_{\lambda_i \in \Lambda_J^c} \left|s_J^{l}(\lambda_i)\right| \ \min_{q_{m}\in \mathcal{P}_{m}^1} \max_{\lambda_i \in \Lambda_J^c} \left| q_{m}(\lambda_i) \right|,
    \end{eqnarray*}
    with $\kappa(\lambda_i) := \|\hat{\mathbf{v}}_i\|_2\|\mathbf{v}_i\|_2$, where $\hat{\mathbf{v}}_i$ and $\mathbf{v}_i \in \mathbb{C}^{n}$ are the left and right eigenvectors associated to $\lambda_i\in\sigma(\mathbf{A})$, with $\hat{\mathbf{v}}_i^{*}\mathbf{v}_i=1$ (see~\cite[Chapter 52]{TrefethenEmbree2005SPB}).
\end{theorem}

As explained in~\cite{DoleanMarchandEtAl2025CAG}, if all the HR values of $N_J^{(l)}$ are close to the eigenvalues in $\Lambda_J$ at iteration $l$, then $\max_{\lambda_i \in \Lambda_J^c} |s_J^{l}(\lambda_i)| \approx 1$.
For the following iterations, the convergence bound then depends only on the other eigenvalues in $\Lambda_J^c$.

The non-linear GMRES convergence rate behavior can be interpreted in local steps with the ratio $\|\mathbf{r}_{l+m}\|_2 / \|\mathbf{r}_l\|_2$.

\section{Domain decomposition preconditioning with deflation}
\label{sec:acc}

To deal with large-scale Helmholtz problems such as~\eqref{eq:linsys}, it can be necessary to resort to domain decomposition methods.
Also, the presence of quasiresonances can hinder the GMRES convergence, and deflation is a useful tool to remove such effect~\cite{DoleanMarchandEtAl2025CAG}.
In this section, we present both methods and their combination.

The computational domain \(\Omega\) is decomposed into \(\mathcal{S}\) overlapping subdomains \(\Omega_s \, \left( s=1,\ldots,\mathcal{S} \right) \), with \(n_s\) the local number of inside dofs.

The one-level method used here is the Optimized Restricted Additive Schwarz (ORAS) preconditioner~\cite{DoleanFryEtAl2025AWR}.
It requires restriction matrices \(\mathbf{R}_s \in \mathbb{R}^{n_s \times N}\), which restrict a global vector to the subdomain \(\Omega_s\), and the transpose \(\mathbf{R}_s^T\) that extends a local vector by zero outside \(\Omega_s\).
Define diagonal matrices \(\mathbf{D}_s \in \mathbb{R}^{n_s \times n_s}\) such that\vspace{-0.1cm}
\begin{eqnarray*}\textstyle
    \mathbf{I} = \sum_{s=1}^{\mathcal{S}} \mathbf{R}_s^T \mathbf{D}_s \mathbf{R}_s \in \mathbb{R}^{N \times N}.
\end{eqnarray*}
Then, by noting \(\mathbf{B}_s \in \mathbb{C}^{n_s \times n_s}\) the local matrices of discretized physical problem, with impedance boundary conditions (BC) on the subdomains interfaces, the ORAS preconditioner reads\vspace{-0.35cm}
\begin{eqnarray*}\textstyle
    \mathbf{M}_{\mathrm{oras}}^{-1} := \sum_{s=1}^{\mathcal{S}} \mathbf{R}_s^T \mathbf{D}_s \mathbf{B}_s^{-1} \mathbf{R}_s.
\end{eqnarray*}

Deflation removes some well-chosen eigenvalues of $\mathbf{A}$ by using projection operators, see e.g.~\cite{GarciaRamosKehlEtAl2020PDM}.
We assume that we have an exact or approximate expression of $n_{\mathrm{def}} \ll N$ eigenvectors, which we store in the columns of a matrix $\mathbf{Z} \in \mathbb{C}^{N \times n_{\mathrm{def}}}$.

\begin{definition}{{(Deflation matrices)}}\label{def:def}
    Let $\mathbf{Z} \in \mathbb{C}^{N \times n_{\mathrm{def}}}$ be a full rank matrix such that $\ker(\mathbf{Z}^{*}) \cap \mathrm{range}({\mathbf{AZ}}) = \{0\}$, and let
    \( \mathbf{E} := \mathbf{Z}^{*}\mathbf{A}\mathbf{Z}\in \mathbb{C}^{n_{\mathrm{def}}\times n_{\mathrm{def}}} \)
        and \( \mathbf{Q} := \mathbf{Z} \mathbf{E}^{-1} \mathbf{Z}^{*}. \)
    The \emph{deflation matrices} are defined as\vspace{-0.1cm}
    \begin{equation*}
        \mathbf{P}_\mathrm{def} := \mathbf{I} - \mathbf{A} \mathbf{Q}
        \quad\text{and}\quad
        \mathbf{Q}_\mathrm{def} := \mathbf{I} - \mathbf{Q} \mathbf{A}.
    \end{equation*}
\end{definition}
The deflation matrices are projectors onto the complement of the deflated subspace.
Applying $\mathbf{P}_\mathrm{def}$ or $\mathbf{Q}_\mathrm{def}$ to $\mathbf{A}$ removes the components associated with the deflated eigenvectors.

To combine deflation with domain decomposition preconditioning, we consider here the adapted deflation preconditioner \(\mathbf{P}_{\mathrm{adef}}\) defined as
\(    \mathbf{P}_{\mathrm{adef}} := \mathbf{M}_{\mathrm{oras}}^{-1} \mathbf{P}_\mathrm{def} + \mathbf{Q}.\)

The deflated and preconditioned problem to solve is then\vspace{-0.1cm}
\begin{align}
    \label{eq:deflatedLinearSystem}
    \left|\ 
    \begin{aligned}
      \text{Find ${\mathbf{x}} \in \mathbb{C}^{N}$ such that $\mathbf{A}\mathbf{P}_\mathrm{adef} {\mathbf{u}} = \mathbf{b}$}, \quad \mathbf{x} = \mathbf{P}_\mathrm{adef} {\mathbf{u}}.
    \end{aligned}
    \right.
\end{align}
To ensure the non-singularity of~\eqref{eq:deflatedLinearSystem}, the required conditions~\cite[Theorem 3.10]{GarciaRamosKehlEtAl2020PDM} are
\begin{equation}\label{eq:singularSystem}
    \ker(\mathbf{Z}^{*}) \cap \mathrm{range}\left(\mathbf{AZ}\right) = \{0\} \quad \text{and} \quad \ker(\mathbf{Z}^{*}) \cap \mathrm{range}\left(\mathbf{M}_{\mathrm{oras}}\mathbf{Z}\right) = \{0\}.
\end{equation}
If the columns of $\mathbf{Z}$ are eigenvectors of $\mathbf{A}$, then the hypothesis of Definition~\ref{def:def} holds and \( \mathbf{E} \) is invertible.
In practice, we assume that the second condition of~\eqref{eq:singularSystem} is satisfied.

\section{Numerical investigations}
\label{sec:num}

\begin{figure}\centering
    \includegraphics[width=0.42\textwidth]{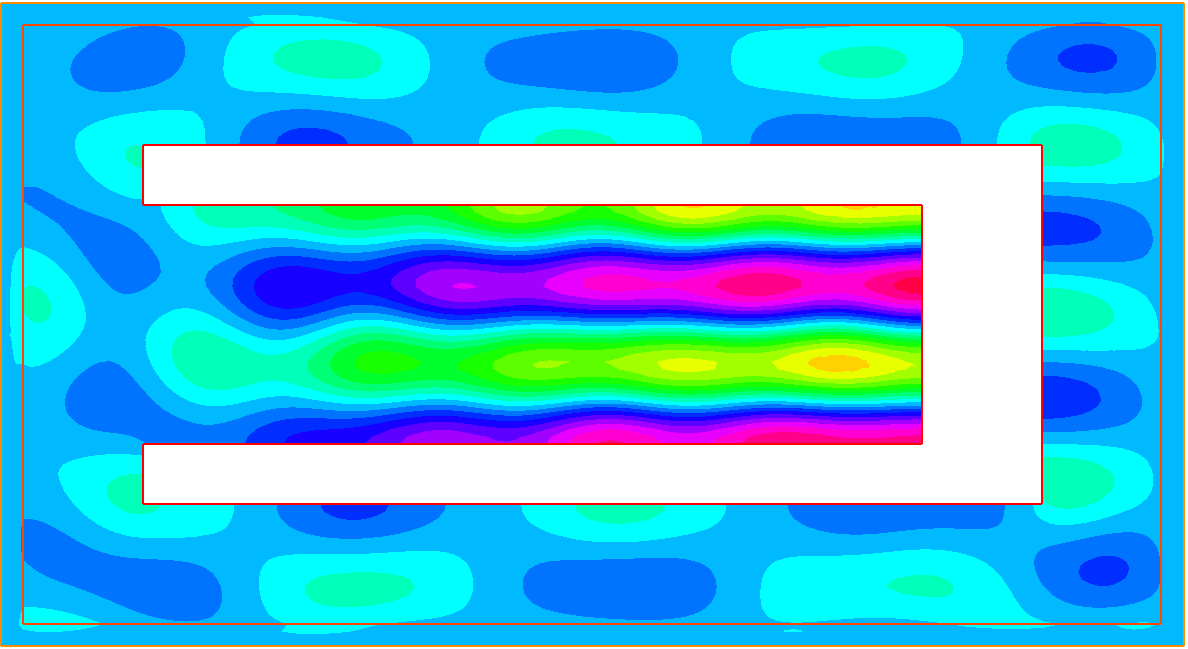}
    \includegraphics[width=0.28\textwidth]{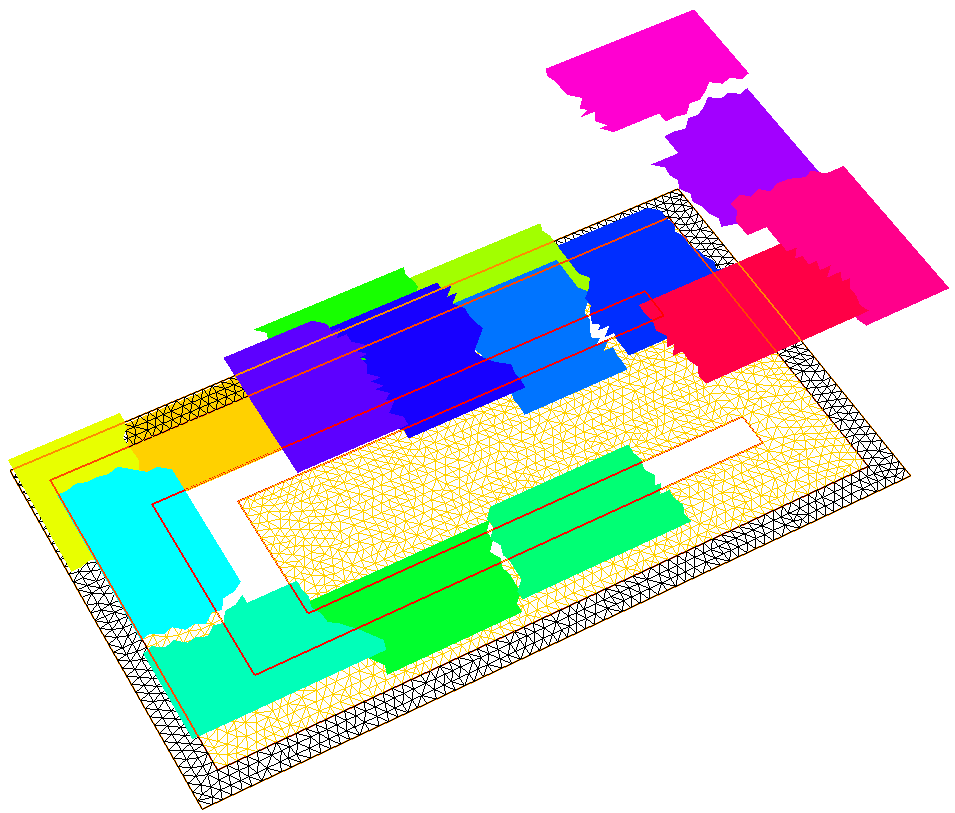}
    \caption{Left: real part of the solution of~\eqref{eq:helmholtzPB}. Right: domain decomposition $(\mathcal{S}=16)$ with the mesh.}
    \label{fig:scatteringPML}
\end{figure}

We consider the scattering of a time-harmonic plane wave with incident angle $\theta$ $u_\mathrm{inc}(x,y) = e^{ik(\cos(\theta) x + \sin(\theta) y)}$, by an obstacle $\mathcal{O}$ in $\mathbb{R}^2$ with Neumann BC on its surface \(\Gamma_\mathrm{obs} \).

The computational area is restricted to a rectangular domain $\Omega_\mathrm{dom} = (-L_x, L_x) \times (-L_y, L_y) \setminus \mathcal{O}$ and surrounded with perfectly matched layers (PMLs) $\Omega_\mathrm{pml}$, of thickness $L_\mathrm{pml}$~\cite{BermudezHervellaNietoEtAl2004EBP}.
The external boundary is denoted by \(\Gamma_\mathrm{ext}\).
The PML formulation consists in a change of variables, with complex-valued stretching functions
\[
    \gamma_x(x) := 1 + {i\sigma_x(x)}/{k}, \quad \gamma_y(y) := 1 + {i\sigma_y(y)}/{k},
\]
where
\(
    \sigma_x(x) := {\mathbbm{1}_\mathrm{pml}(x)}/\left( {L_\mathrm{pml} - |x| + L_x} \right)\) and \(
    \sigma_y(y) := {\mathbbm{1}_\mathrm{pml}(y)}/\left( {L_\mathrm{pml} - |y| + L_y} \right)
\).
The scattered field $u$ then satisfies:
\begin{align}\label{eq:helmholtzPB}
    \left\{
    \begin{aligned}
        {\partial_x} \left( {\gamma_y}/{\gamma_x} {\partial_x}u \right) +
        {\partial_y} \left( {\gamma_x}/{\gamma_y} {\partial_y}u \right) +
        \gamma_x \gamma_y k^2 u & = 0 \quad & \text{in } \Omega_\mathrm{dom} \cup \Omega_\mathrm{pml}, \\
        \partial_{\boldsymbol{n}} u & = -\partial_{\boldsymbol{n}} u_\mathrm{inc} \quad & \text{on } \Gamma_\mathrm{obs}, \\
        u & = 0 \quad & \text{on } \Gamma_\mathrm{ext}.
    \end{aligned}
    \right.
\end{align}

The obstacle $\mathcal{O}$ is an open rectangular cavity (see Fig.~\ref{fig:scatteringPML}) of length $L_\mathcal{O}$ and opening width $l_\mathcal{O}$.
For this configuration, quasimodes are expected (see Definition~\ref{def:quasimodes}).

They induce small eigenvalues~\cite{DoleanMarchandEtAl2025CAG} when the wave numbers are close to the resonance frequencies of a closed rectangular cavity with Neumann BC on all sides except the left edge (Dirichlet), given by
\begin{equation}\label{eq:quasimodes}
    k_{n,m} = \pi \sqrt{{{(m + 1/2)}^2}/{L_\mathcal{O}^2} + {n^2}/{l_\mathcal{O}^2}}, \quad m,n > 0.
\end{equation}
The eigenvectors associated with the small eigenvalues resulting from the quasimodes are close to the eigenmodes of the corresponding closed configuration.

Discretization by finite elements leads to a linear system $\mathbf{A} \mathbf{u} = \mathbf{b}$, where the matrix $\mathbf{A}$ is complex, has complex-valued eigenvalues, and is non-normal. 

The following parameters are chosen: 
$L_\mathcal{O} = 1.3$, $l_\mathcal{O} = 0.4$, $\theta = {4\pi}/{10}$, \text{$P_3$ elements}, $\mathtt{tol} = 10^{-6}$.
For the first set of experiments, the number of dofs per wavelength is fixed to $10$, for a total of about $N \approx 3.10^4$.
The second experiments are performed with $5$ dofs per wavelength and $N \approx 1.10^5$, at a higher wave number.

For running the numerical experiments, we used the framework \texttt{FFDDM}~\cite{FFD:Tournier:2019} of \texttt{FreeFEM}~\cite{MR3043640}.
The first experiments illustrate how quasimodes affect GMRES convergence with ORAS preconditioning.
The second set of experiments investigates the influence of quasimodes and their deflation when coarse space preconditioners are used.

\subsubsection*{Influence of quasimodes on the GMRES convergence with ORAS preconditioning}

\begin{figure}\centering
    \includegraphics[width=\textwidth]{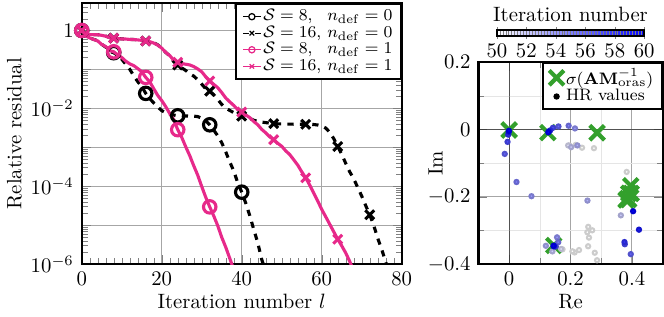}
    \caption{GMRES convergence with ORAS preconditioning for different number of subdomains \(\mathcal{S}\) (left). Spectrum of $\mathbf{A}\mathbf{M}_{\mathrm{oras}}^{-1}$ and HR values trajectories for \(\mathcal{S}=16\) (right).}
    \label{fig:OrasDefSpectrum}
\end{figure}

We consider the wave number \(k=23.591\) to be close to \(k_{0,3}\) in~\eqref{eq:quasimodes}.
We investigate how the quasimodes affects the GMRES convergence $\mathcal{S}=8,16$ subdomains in Fig.~\ref{fig:OrasDefSpectrum} (right).

With \(\mathcal{S}=8\) and \(16\) (black), a stagnation phase is observed during the mid-iterations, before a fast convergence phase.
This behavior is typical for quasiresonant Helmholtz problems~\cite{DoleanMarchandEtAl2025CAG}, and Theorem~\ref{thm:harmonicRitzValues} gives an interpretation with HR values.

The spectrum of \(\mathbf{A}\mathbf{M}_{\mathrm{oras}}^{-1}\) for \(\mathcal{S}=16\) (Fig.~\ref{fig:OrasDefSpectrum}, right) shows that there are a few small eigenvalues with a small imaginary part, which are associated with quasimodes.
The positions of the HR values around iteration $50$ are plotted in $\bullet$ gray, and they become increasingly blue as the iteration increases up to $60$.

We observe that for $\mathcal{S}=16$, the plateau of the residual matches the iterations where the HR values approach the two smallest eigenvalues related to quasimodes.
At iteration $60$, these small eigenvalues are well approached by HR values.
From this iteration, Theorem~\ref{thm:harmonicRitzValues} indicates that the influence of the small eigenvalues on the convergence is removed, so the convergence rate increases.

This analysis reveals that the plateau observed in the middle of the convergence is due to the quasimodes and not only to the number of subdomains.
However, it should be noted that when the cavity is divided into a greater number of subdomains, the plateau lasts longer.

In Fig.~\ref{fig:OrasDefSpectrum} (left), the same experiments were conducted by deflating the smallest eigenvalue, using an approximation of the closed cavity mode, extended by zero outside the cavity (pink).
These experiments confirm that the plateaus are due to quasimodes and disappear when the associated modes are deflated.

\subsubsection*{Influence of quasimodes with coarse spaces preconditioners}

We now consider two-level preconditioners with coarse spaces (CS) built with Dirichlet-to-Neumann (DtN) eigenproblems~\cite{DoleanFryEtAl2025AWR}, and we investigate the influence of quasimodes on the GMRES convergence.
The same numerical experiments have been carried out with H-GenEO CS~\cite{DoleanFryEtAl2025AWR} and have yielded similar results (not shown here for sake of brevity).

\begin{figure}
    \includegraphics[width=\textwidth]{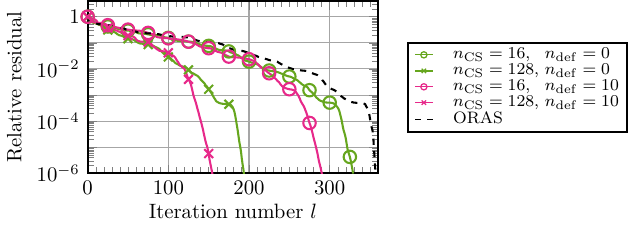}
    \caption{GMRES convergence with DtN CS, with and without additional deflation contribution for $k=102.11$ and $\mathcal{S}=16$.}\label{fig:DtNGenEO}
\end{figure}

In Fig.~\ref{fig:DtNGenEO}, the evolution of the relative residual is plotted for $k=102.11$ (close to $k_{0,13}$) with $\mathcal{S}=16$ and two coarse space sizes ($n_{\mathrm{CS}}$).

In the residual history (green), the first stagnation phase is related to the number of subdomains, and it is significantly reduced by the coarse space.
However, several late plateaus remain.
By deflating the small eigenvalues associated with the quasimodes, the late stagnation phases disappear (pink).
So, these late plateaus are due to the quasimodes.
The GMRES convergence preconditioned by this type of coarse space is therefore also hindered by the quasimodes.

\begin{center}\vspace{-0.75cm}
\small
\begin{table}[h!]\centering
\begin{tabular}{|c||c|c|c||c|c|c|}
    \hline
     & \multicolumn{3}{c||}{$n_{\mathrm{CS}} + n_{\mathrm{def}} = 32$} & \multicolumn{3}{c|}{$n_{\mathrm{CS}} + n_{\mathrm{def}} = 144$} \\
    \hline
    \shortstack{CS \\ composition} 
    & \shortstack{$n_{\mathrm{CS}}=32$\\ $n_{\mathrm{def}}=0\phantom{0}$}
    & \shortstack{$n_{\mathrm{CS}}=0$\\ $n_{\mathrm{def}}=32$}
    & \shortstack{$n_{\mathrm{CS}}=16$\\ $n_{\mathrm{def}}=16$}
    & \shortstack{$n_{\mathrm{CS}}=144$\\ $n_{\mathrm{def}}=0\phantom{00}$}
    & \shortstack{$n_{\mathrm{CS}}=0\phantom{00}$\\ $n_{\mathrm{def}}=144$}
    & \shortstack{$n_{\mathrm{CS}}=128$\\ $n_{\mathrm{def}}=16\phantom{0}$} \\
    \hline
    $\#\mathrm{iter}$ & $302$ & $311$ & $273$ & $180$ & $236$ & $141$ \\
    \hline
\end{tabular}
\caption{Number of GMRES iterations for two second-level sizes comparing the contributions of the DtN CS, deflation, and the combination of the two, for $k=102.11$, $\mathcal{S}=16$.}
\label{tab:DtNGenEO}
\end{table}
\end{center}
\vspace{-1.16cm}

To deal with large-scale problems with quasimodes, it can be interesting to combine standard coarse space preconditioning with deflation of quasimodes.
Tab.~\ref{tab:DtNGenEO} shows the number of GMRES iterations required to reach \(\mathtt{tol}=10^{-6}\), with two fixed coarse space sizes ($n_{\mathrm{CS}} + n_{\mathrm{def}} = 32$ or $144$), by varying the contributions of the DtN CS and of the deflation.
The deflation only is the least efficient strategy, because the local deflated vectors cannot handle the global exchange of information ensured by the coarse spaces (first stagnation phase in Fig.~\ref{fig:DtNGenEO}).
Using only the DtN CS is efficient, but as previously mentioned, it is sensitive to the quasimodes.
The third and sixth columns correspond to using both DtN CS and deflation.
This combined approach is the most efficient.

For geometries with explicit or approximate expressions of the quasimodes, adding a few deflated vectors to the coarse spaces can thus be worthwhile.
Besides, in this case, it is a priori less costly numerically than increasing the size of a coarse space.

\section{Conclusion}

This work highlighted the impact of quasimodes on the GMRES convergence in the context of domain decomposition methods.

The study with HR values showed that dividing a quasiresonant cavity into several subdomains generates small eigenvalues that hinder GMRES convergence.
The more the cavity is divided, the longer the stagnation phase.

Preconditioning with DtN coarse spaces is not sufficient to overcome the impact of quasimodes.
A combination of these preconditioners and deflation of quasimodes is robust against the number of subdomains (coarse spaces contribution), and the plateaus due to the small eigenvalues related to quasimodes (deflation contribution).

\ethics{Acknowledgements}{
This work was funded in part by the \textit{ANR JCJC project WavesDG} (research grant ANR-21-CE46-0010) and by the \textit{Agence de l'Innovation de Défense} [AID] through \textit{Centre Interdisciplinaire d'Etudes pour la D\'efense et la S\'ecurit\'e} [CIEDS] (project 2022 ElectroMath).\newline
The authors would like to thank Pierre-Henri Tournier for his help with FreeFEM implementation.
}

\vspace{-0.3cm}
\bibliographystyle{spmpsci}
\bibliography{bibliography}
\end{document}